\documentclass[12pt]{amsart}

\usepackage{graphicx}
\usepackage{amsmath,amssymb,amsfonts}
\usepackage{amsthm}
\usepackage{amssymb}
\usepackage{mathrsfs}
\usepackage{xcolor}%
\usepackage{textcomp}%
\newtheorem{theorem}{Theorem}
\newtheorem{proposition}[theorem]{Proposition}
\newtheorem{remark}{Remark}

\usepackage[utf8]{inputenc}
\usepackage[dvipsnames]{xcolor}
\usepackage{lmodern}
\usepackage{color}
\usepackage{hyperref}
\usepackage{letltxmacro}
\hypersetup{
    colorlinks=true, 
    linktoc=all, 
    linkcolor=blue, 
   citecolor=blue,
    filecolor=blue,
    urlcolor=blue
}

\newcommand{\SL}{\operatorname{SL}}
\newcommand{\GL}{\operatorname{GL}}
\newcommand{\SO}{\operatorname{SO}}

\newtheorem{lemma}[theorem]{Lemma}

\newtheorem{acknowledgement}[theorem]{Acknowledgement}
\def\R{\mathbb R}

\title{Conjugacy limits of certain subgroups in $\SL(2,\mathbb{R})\ltimes\mathbb{R}^2$}
\author{Manoj Choudhuri}
\address{Institute of Infrastructure, Technology, Research And Management, Near Khokhara Circle, Maninagar (East), Ahmedabad 380026, Gujarat, India.}
\email{manojchoudhuri@iitram.ac.in}
\author{C. R. E. Raja}
\address{Theoretical Statistics and Mathematics Unit, Indian Statistical Institute, 8th Mile, Mysore Road, R. V. College Post, Bangalore 560059, Karnatake, India.}
\email{creraja@isibang.ac.in}

\begin{document}  

\begin{abstract} 
We study conjugacy limits of certain of subgroups inside $\SL(2,\R)\ltimes\R^2$. These subgroups have a common feature that any two in the same category are conjugates of each other.
\end{abstract}

\subjclass[2020]{Primary: $22$E$15$, $54$H$11$ Secondary: $37$B$02$}
 \maketitle
\tableofcontents

\section{Introduction}
Let $G$ be a locally compact group and  $\text{Sub}_G$ be the space of closed subgroups of $G$ equipped with the Chabauty topology. Chabauty topology on $\text{Sub}_G$ was introduced by Claude Chabauty in $1950$ (\cite{CC}). Chabauty's purpose in introducing this topology was to generalise Mahler's compactness criterion for lattices in $\R^n$ to a larger class of locally compact groups. The Chabauty topology on $\text{Sub}_G$
is generated by the subsets $\{ H \in \text{Sub}_G \mid H\cap K = \emptyset \}$, $\{ H\in \text{Sub}_G \mid H \cap U \not = \emptyset \}$, 
where $K \subset G$ is compact and $U \subset G$ is open, and $\text{Sub}_G$ becomes a compact metrizable space with respect to this topology. A sequence $H_n\in\text{Sub}_G$ converges to $H\in \text{Sub}_G$ if and only if for any $h\in H$ there is a sequence $h_n\in H_n$ such that $h_n$ converges to $h$ and for any sequence $x_{k_n}\in H_{k_n}$, with $k_{n+1}>k_n$ and $x_{k_n}\to x$, we have $x\in H$.

Since its introduction, the Chabauty topology has drawn attention of many mathematicians over the years. Besides its wide range of applicability such as in the study of discrete subgroups of semi-simple Lie groups, hyperbolic manifolds, etc., a lot of attentions has been paid just to understand the structure of the space $\text{Sub}_G$, and it is known completely only for a few groups. For example, it is know from the work of  Hubbard and Pourezza (\cite{PH}) that $\text{Sub}_\mathbb C$ is homeomorphic to the four-dimensional sphere $S^4$. In a more recent work, the structure of $\text{Sub}_G$ for certain low-dimensional Lie groups with a special attention to the $3$-dimensional Heisenberg group have been studied by Bridson, Harpe and Kleptsyn in \cite{BHK}. 

In order to understand the space $\text{Sub}_G$, one has to understand certain subspaces of it formed by various classes of closed subgroups and their closures. So, studying limits of subgroups inside $\text{Sub}_G$ becomes essential. Studying limits of subgroups has other usefullness as well, especially, in the study of geometric transition in which one studies the deformation of one geometry to another using the conjugacy limits of subgroups (cf. \cite{CDW}). One may look at \cite{CDW} for various notions of limits in the space $\text{Sub}_G$ as well. The study of conjugacy limits has been an active area of research in order to understand the space $\text{Sub}_G$ and for its connections with geometry. For example, see a series of papers by Arielle Leitner: \cite{L3}, \cite{Ln}, \cite{L3L}. In \cite{L3} and \cite{Ln}, the conjugacy limits of diagonal subgroups are discussed, whereas, \cite{L3L} discusses the cunjugacy limits of connected subgroups. See also the references cited there in these papers for more related works and connections with geometry.

In this article, we study the conjugacy limits of various types of subgroups inside $G=\SL(2,\R)\ltimes \R^2$. Note that $G$ is neither semi-simple nor solvable. The type of subgroups we consider include Levi subgroups, maximal compact subgroups, maximal diagonalizable subgroups, Borel subgroups and unipotent subgroups. Recall that these subgroups have a special feature that any two member of the same type are conjugate of each other. So, studying conjugacy limits of any of these subgroups is equivalent to studying the closure of the orbits of a subgroup under the conjugation action of $G$ on $\text{Sub}_G$. We got interested to find these orbit closures while investigating about the existence of invariant random subgroups supported on these orbits (\cite{CR1}). Theorem $1.1$ of \cite{CR1} shows that the set of Levi subgroups of $\SL(2,\R)\ltimes \R^2$ is not closed, and this was our starting point of investigation about the conjugacy limits of the class of subgroups considered in this article. 

 Given a closed subgroup $H$ of $G$, by a conjugacy limit of $H$ in $\text{Sub}_G$, we mean the limit of $(g_nHg_n^{-1})$ for a sequence $(g_n)$ in $G$. Note that $(g_nHg_n^{-1})$ always has a convergent subsequence as $\text{Sub}_G$ is compact. Let us denote by $L$ the subgroup $\SL(2,\R)$ of $G=\SL(2,\R)\ltimes \R^2$. Then $L$ is a Levi subgroup of $G$. The following theorem gives a complete description of all the limit points of the set of all Levi subgroups of $G$ which we denote by $L_G$. We write any element of $G=\SL(2,\R)\ltimes \R^2$ as $(g,v)$ with $g\in\SL(2,\R)$ and $v\in\R^2$. Also, for any two elements $(g,v)$ and $(g',v')$, the group operation is defined as follows:
\[
(g,v)(g',v')=(gg',v+g(v')).
\]
We identify $\SL(2,\R)$ and $\R ^2$ as closed subgroups of $G$ by the maps $g \mapsto (g, (0,0))$ 
and $v\mapsto (I,v)$ where $I$ denotes the $2\times 2$ identity matrix.  For any $c\in \R$, we 
define a unipotent subgroup $N^+ _c$ by 
\[
N^+ _c=\left\{ \big(\begin{pmatrix} 1&x\\ 0&1 \end{pmatrix}, (cx, 0) \big):x\in\R\right\}.
\]
For $c=0$, we denote $N^+_c$ by $N^+$, which is identified (and isomorphic) to 
the group of unipotent upper triangular matrices in $\SL(2,\R)$.  
\begin{theorem}\label{limit-L}
Any subgroup in $\overline {L_G}$ is either a Levi subgroup of $G$ or a conjugate of 
$N^+\ltimes \R^2$. 
\end{theorem}
Let $K$ be a maximal compact subgroup of $G$; for example, $K=\SO(2,\R)$, and let $K_G$ be the set of all maximal compact subgroups of $G$, i.e., \[ K_G= \{ gKg^{-1} \mid g\in G \} .\] Also, consider the following subgroups of $G$:
\[
V_c=\{ (I, (t, ct)) \mid t \in \R \},\ c\in\R, ~~
V_\infty=\{ (I, (0, t)) \mid t \in \R \}
\]
and 
\[ \tilde{N}^+=\left\{ \left(\pm \begin{pmatrix}
    1&s\\ 0&1
\end{pmatrix},(0,0)\right):s\in\R\right\}.
\]
It may be noted that $V_c , c\in \R \cup \{\infty \}$ is the set of all one-dimensional subspaces of $\R ^2$, the projective space $\mathcal{P} ^1$.  
\begin{theorem}\label{limit-K}
The closure of $K_G$ consists of conjugates of the subgroups $V_c,c\in \R\cup\{\infty\}$ and $\tilde{N}^+$, in addition to the conjugates of $K$.  
\end{theorem}
The next result identifies limit points of the set of maximal diagonalizable subgroups of $\SL(2,\R)\ltimes \R^2$. Let $A= \left\{ \begin{pmatrix} a &0 \\ 0& \frac{1}{a} \end{pmatrix} \mid a \in \R^+ \right\}$ and $D_G$ be the set of all maximal diagonalizable subgroups of $G$, i.e., \[D_G = \{ gA g^{-1} \mid 
g\in SL(2, \R)\ltimes \R ^2 \}.\]  We now introduce the following one-dimensional subgroups
\[
V_{(a,b,c)}= \left\{ \left( \begin{pmatrix} 1& at \\ 0&1 \end{pmatrix} ,\left( ct+\frac{abt^2}{2}, bt\right) \right ) \mid t\in \R \right\} 
\]
for any $a, b, c \in \R$ and not all of them are zero. It may be noted that $V_{a,b,c}$'s are the one-parameter subgroups in the $3$-dimensional Heisenberg group and $V_{(1,0,c)}= N^+_c$, $V_{(0,1,c)}= V_c$ and $V_{(0,1,0)} =V_\infty$.  
\begin{theorem}\label{limit-D}
The closure of $D_G$ consists of conjugates of \[V_{(a, b, c)},\ a, b, c \in \R\] with $a =1$ or $b=1$, in addition to the conjugates of the diagonal group $A$.
\end{theorem}
Let $B$ denote the group of all upper triangular matrices in $\SL(2,\R)$. Then $B$ is a Borel subgroup of $\SL(2,\R)$ and $B\ltimes \R^2$ is a Borel subgroup of $G=\SL(2,\R)\ltimes \R^2$. Then, using the Iwasawa decomposition in $\SL(2,\R)$, it is easy to see that, for conjugacy limits of $B\ltimes \R^2$, we only need to consider conjugation by elements of $K=SO(2,\R)$, hence conjugates of $B$ form a closed subset. We, therefore, consider conjugacy limits of $B$ inside $G$ and identify all possible conjugacy limits of $B$.
\begin{theorem}\label{limits-B}
The closure of $\{ gBg^{-1} \mid g\in SL(2, \R) \ltimes \R ^2 \}$ consists of conjugates of $B$, $N^+\ltimes \{ (t,0) \mid t\in \R \} ~(\simeq \R ^2)$ and $\R^2$. 
\end{theorem}

The next result is concerned about the conjugacy limits of the unipotent subgroup $N^+$. 
\begin{theorem}\label{limits-U}
The closure of $\{ gN^+g^{-1} \mid g\in SL(2, \R) \ltimes \R ^2 \}$ consists of conjugates of $N^+$ and $V_0$.
\end{theorem}

\section{Preliminary}
In this section, we prove certain technical results which are needed for what follows.

\begin{proposition}\label{cs}
Let $G$ be an exponential Lie group and $H_n \to H$ in $\text{Sub}_G$.  Let $\mathcal{ G}$, $\mathcal{H} _n$ and $\mathcal{H}$ be Lie algebras of $G$, $H_n$ and $H$ respectively.  If $H_n$'s are connected, then $H$ is also connected, $\mathcal{H}_n \to \mathcal{H}$ in Sub$_\mathcal{G}$ and $\limsup ~\dim (H_n) = \dim (H)$.   

\end{proposition}

\begin{remark}
 It may be noted that $H$ need not always be connected even when $H_n$ are all connected and have the same dimension: e.g., last case of Section $4$ provides a sequence of conjugate subgroups whose limit point has two connected components. 
\end{remark}
\begin{proof}
Let $V$ be a limit of $\mathcal{H}_n$ in Sub$_\mathcal{G}$.  By passing to a subsequence, we may assume that $\mathcal{H}_n \to V$ in Sub$_\mathcal{G}$. Then $V$ is a vector subspace. It is easy to see that $V$ is a Lie subalgebra of $\mathcal{G}$ as well.  For $v\in V$, there are $v_n \in \mathcal{H} _n$ such that 
$v_n \to v$, hence $tv_n \to tv$ for all $t\in \R$.  Let $\exp \colon \mathcal{G} \to G$ be the exponential map.  
Then $\exp (tv_n) \to \exp (tv)$ for all $t\in \R$.  Since $H_n$ converges to $H$, $\exp (tv) \in H$ for all $t\in \R$.  This implies 
that $v\in \mathcal{H}$.  Thus, $V \subset \mathcal{H}$.

Assume that $\exp \colon \mathcal{G} \to G$ is a diffeomorphism and $H_n$ are connected.  Then $\exp \colon \mathcal{H} _n \to H$ is also a diffeomorphism for all $n$.  
Since $H_n \to H$, for $h\in H$, there are $h _n \in H_n$ such that $h _n \to h$. Since $G$ and $H_n$ are exponential, there are $X_n \in \mathcal{H} _n$ and $X\in \mathcal{G}$ be such that $\exp (X_n) = h_n$ and $\exp (X) =h$, hence $\exp (X_n) \to \exp (X)$.  Since $\exp$ is a diffeomorphism, we get that 
$X_n \to X$.  Since $\mathcal{H} _n \to V \subset \mathcal{H}$, we get that $X\in \mathcal{H}$.  Thus, $\exp (\mathcal{H} ) =H$.  Therefore, $H$ is connected.  

If $Y\in \mathcal{H}$, then $\exp (Y) \in H$.  Since $H_n \to H$, there are $g_n \in H_n$ such that $g_n \to \exp (Y)$.  Since all $H_n$ are exponential, there are $Y_n \in \mathcal {H}_n$ such that $\exp (Y_n) = g_n$.  Thus, $\exp (Y_n) \to \exp (Y)$.  Since $\exp$ is a diffeomorphism, we get that $Y_n \to Y$.  Since $\mathcal{H} _n \to V$, $\mathcal{H}\subset V$ 
but $V\subset \mathcal{H}$, hence $V= \mathcal{H}$.  Thus, $\limsup ~\dim (H_n) = \dim (H)$.    
\end{proof} 

The next lemma is a key ingredient in the proof of the main results.

\begin{lemma}\label{el1} 
Let $(g_n)$ and $(v_n)$ be sequences in $\GL (n,\R)$ and $\R^n$ respectively such that  
\[(I,v_n)\ (g_n,0)\ (I,v_n)^{-1}=(g_n,v_n-g_n(v_n))\rightarrow (g,v) \] 
for some $g\in\GL(n,\R)$ and $v\in \R ^n$. Then, either $1$ is an eigenvalue of $g$ or $(I-g)$ is invertible and $v_n \to (I-g)^{-1}(v)$.  
In particular, if $g\in SL(2, \R)$, then $g$ is unipotent or $v_n$ converges. 
\end{lemma}

\begin{proof} 
Suppose $(g_n)$ and $(v_n)$ are sequences in $\GL (n,\R)$ and $\R ^n$ respectively such that 
\[(I,v_n)\ (g_n,0)\ (I,v_n)^{-1}=(g_n,v_n-g_n(v_n))\rightarrow (g,v)\] for some $g\in\GL(n,\R)$ and $v\in \R^n$.  
Then $g_n \to g $ in $GL(n, \R )$ and  $(I-g_n)(v_n)\rightarrow v$ as $n\rightarrow\infty$.  
If $1$ is not an eigenvalue of  $g$, then $I-g$ is invertible.  Since $g_n \to g$, $I-g _n \to I-g$.  Since invertible matrices form an open set, $I-g_n$ is invertible for all large $n$ and hence $(I-g_n)^{-1} \to (I-g)^{-1}$.  As $(I-g_n)(v_n)\rightarrow v$, we get that 
$v_n  \to (I-g)^{-1} (v)$.  

If $g\in SL( 2, \R )$ and $v_n$ does not converge, then one is an eigenvalue of $g$.  Since $g\in SL(2, \R)$, one is the only eigenvalue of $g$, hence $g$ is an unipotent transformation.

\end{proof}

\section{ Conjugacy limits of Levi subgroups}

Let $G =\SL(2,\R)\ltimes \R^2$, $L = \SL(2,\R)$ and $L_G$ be the set of all semisimple Levi subgroups inside $G$. Then $L_G$ consists of conjugates of $L$ by elements of $G$ (cf. Theorem 3.18.3 of \cite{Va}). We denote elements in $G$ by $(g,v)$ with $g\in\SL(2,\R)$ and $v\in\R^2$. Now, if $(G_n)$ is a sequence of Levi subgroups in $L_G$, then each $G_n$ is given by a conjugate of $L$ by some element of $G$, i.e.,
\[
G_n=\ (g_n,v_n)\ L\ (g_n,v_n)^{-1},
\]
where $(g_n,v_n)^{-1}$ denotes the inverse of $(g_n,v_n)$ in $G$.

If we conjugate $L$ by elements of the form $(g,(0,0))$, the resultant is $L$ itself. Also, note that $(g,v)=(I,v)\ (g,(0,0))$, where $I$ is the identity element in $\SL(2,\R)$. Therefore, to determine all the limit points of $L_G$, it is enough to consider the conjugacy limits of the form
\[
(I,v_n)\ L\ (I,v_n)^{-1}.
\]

\begin{proof}[Proof of Theorem \ref{limit-L}]
Assume that $(I,v_n)\ L\ (I,v_n)^{-1}$ converges to $H$ in $\text{Sub}_G$; we now determine $H$. The proof is divided into following caes. 

\noindent {\bf Case $1$.} Suppose $v_n\rightarrow v$ as $n\rightarrow\infty$. In this case, it is easy to see that $(I,v_n)\ L\ (I,v_n)^{-1}$ converges to the closed subgroup
$(I,v)\ L\ (I,v)^{-1}$, and hence the limit $H$ is a Levi subgroup.

\smallskip

We now assume that $(v_n)$ is unbounded. 
\\ 
\noindent {\bf Case $2$.} Suppose $v_n=\alpha_n\ e_1$, where $e_1$ is the unit vector on $x$-axis and $|\alpha_n|\rightarrow\infty$. Now, given any $s\in\R$, if $g_n=\begin{pmatrix} \frac{\alpha_n-s}{\alpha_n} & 0\\ 0 & \frac{\alpha_n}{\alpha_n-s} \end{pmatrix}$, then 
\[
(I,v_n)\ (g_n,(0,0))\ (I,v_n)^{-1}=(g_n,v_n-g_n(v_n))=(g_n,(s,0))\rightarrow (I,(s,0))
\]
as $n\rightarrow\infty$. Again, taking $g_n=\begin{pmatrix} 1 & 0\\\frac{-s}{\alpha_n}  & 1 \end{pmatrix}$, we see that 
\[
(I,v_n)\ (g_n,(0,0))\ (I,v_n)^{-1}=(g_n,v_n-g_n(v_n))=(g_n,(0,s))\rightarrow (I,(0,s))
\]
as $n\rightarrow\infty$. So, $\R^2$ is contained in the limit subgroup $H$. Also, any element of the unipotent subgroup $N^+$ fixes $v_n$. Therefore, $N^+\subset H\cap SL(2, \R)$.  It follows from Lemma \ref{el1} that if $g\in SL(2, \R ) \cap H$, then $g$ is unipotent. Thus, $SL(2, \R) \cap H$ is a subgroup of $G$ consisting of unipotent transformations and $SL(2, \R) \cap H$ contains $N^+$ which is one of the maximal group consisting of unipotent transformations, hence $SL(2, \R) \cap H = N^+$.  Thus, $H= N^+\ltimes \R ^2$.  

\noindent {\bf Case $3$.} Let $v_n=(a_n,b_n)$ with $||v_n||$ unbounded.  Let $\alpha _n = ||v_n||$.  Then there are $A_n \in SO(2, \R)$ such that 
\[
v_n=\alpha_n A_n(e_1)=A_n(\alpha_n\ e_1) = A_n (w_n) 
\] where $w_n= \alpha _n e_1$.  
By passing to a subsequence, we may assume that $A_n$ converges to $A$ in $\SL(2,\R)$.  Now,
\begin{align*}
(I,v_n)L(I,-v_n)&=(A_n,(0,0))(I,w_n)(A_n^{-1},(0,0))L(A_n,(0,0))(I,-w_n)(A_n^{-1},(0,0)\\&=(A_n,(0,0))\big((I,w_n)L(I,-w_n)\big)(A_n^{-1},(0,0))\\&
\rightarrow AHA^{-1}
\end{align*}
by Case $1$.   
\end{proof}

\section{Conjugacy limits of maximal compact subgroups}
We now identify the limit points of the set of maximal compact subgroups. In other words, we find the closure of the orbit of $K=\SO(2,\R)$ under the conjugacy action of $G$ on $\text{Sub}_G$. 
\begin{proof}[Proof of Theorem \ref{limit-K}]
Let us recall that $K$ is given by  
\[
K:=\ \left\{ \begin{pmatrix}
    \text{cos}\ \theta & \text{sin}\ \theta\\
    -\text{sin}\ \theta & \text{cos}\ \theta
\end{pmatrix}\ :\ \theta\in\R \right\}. 
\]
In order to find the limit points of 
\[
K_G=\{ (g, v) K (g, v) ^{-1}\mid g \in SL (2, \R), v\in \R ^2 \},
\]
it is sufficient to identify the limit points up to conjugation. Using Iwasawa decomposition, i.e.,  $SL(2, \R) = N^+AK$, we may assume that $g\in N^+A$.  Let $H$ be a limit point of $K_G$. Then there exist sequences $$\mathfrak{a}_n = \begin{pmatrix} a_n & 0 \\ 0 & a_n^{-1} \end{pmatrix} \in A, ~ \mathfrak{s}_n  = \begin{pmatrix} 1 & s_n \\ 0 &1 \end{pmatrix} \in N^+ $$ and 
$v_n = (\alpha _n , \beta _n)\in \R ^2$ such that $(\mathfrak{s}_n \mathfrak{a}_n, v_n)K(\mathfrak{s}_n \mathfrak{a}_n, v_n) ^{-1} \to H$. By passing to a subsequence and replacing $\mathfrak{a}_n$ or $\mathfrak{s}_n$ by their conjugates, we may assume that $a_n \to \infty$ or $a_n =1$ and $s_n \to \infty$ or $s_n =0$.  

Let $(g, v)\in H$.  Then there exists $$\Theta _n = \begin{pmatrix}
    \cos \theta_n & \sin \theta_n\\
    -\sin \theta_n & \cos \theta_n
\end{pmatrix} \in K$$ such that $(\mathfrak{s}_n \mathfrak{a}_n, v_n)\Theta _n(\mathfrak{s}_n \mathfrak{a}_n, v_n) ^{-1} \to (g, v)$

Now,
\[
(\mathfrak{s}_n \mathfrak{a}_n, v_n) \Theta _n(\mathfrak{s}_n \mathfrak{a}_n, v_n) ^{-1} = (T_n , (x_n, y_n)),
\]
where 
\[
T_n = \begin{pmatrix} \cos \theta _n -\frac{s_n}{a_n ^2}\sin \theta _n  & (\frac{s_n ^2}{a_n^2}+a_n ^2) \sin \theta _n \\ -\sin \theta _n / a_n ^2 & \cos \theta _n +\frac{s_n}{a_n ^2}\sin \theta _n \end{pmatrix}
\]
and 
\[
x_n = \alpha _n \big(1-\cos \theta _n +{\frac{s_n}{a_n ^2}}\sin \theta _n \big) +\beta _n  \big(\frac{s_n^2}{a_n ^2}+a_n ^2\big) \sin \theta _n,
\]
\[
y_n = -\alpha _n \frac{\sin \theta _n}{a_n ^2} + \beta _n \big(1-\cos \theta _n -\frac{s_n}{a_n ^2}\sin \theta _n \big).
\]

We now assume that $v_n \to \infty$. As $\big(T_n , (x_n , y_n)\big) \to (g, v)$, it follows from Lemma \ref{el1} that $g$ is unipotent.  

\noindent {\bf Case $1$:} We now consider the case when ${\frac{\alpha _n}{\beta _n}} \to 0$.  Since $v_n$ diverges, $\beta _n \to \infty$.  Let $x,y\in \R$ be such that $v= (x,y)$.  Then
\[
x_n = \alpha _n \big(1-\cos \theta _n +\frac{s_n}{a_n ^2}\sin \theta _n \big) +\beta _n \big(\frac{s_n^2}{a_n ^2}+a_n ^2\big) \sin \theta _n  \to x.
\]
Since $T_n \to g$, $\cos \theta _n -\frac{s_n}{a_n ^2}\sin \theta _n $ converges.  Therefore, 
\[
x_n = \beta _n \Big[\frac{\alpha _n }{ \beta _n}\big(1-\cos \theta _n +\frac{s_n}{ a_n ^2}\sin \theta _n \big)
+  \big(\frac{s _n^2}{a_n ^2}+a_n ^2\big) \sin \theta _n\Big]
\]
converges implies that $\big(\frac{s_n^2}{a_n ^2}+a_n ^2\big) \sin \theta _n \to 0$. Since $g$ is unipotent, it follows that $g= I$. 

Also, for any $t\in \R$, we can choose $\zeta _n \to 0$ so that $ \beta _n (\frac{s_n^2}{a_n ^2}+a_n ^2) \sin \zeta _n  \to t$.  Then $\beta _n\frac{s_n }{a_n ^2}\sin \zeta _n \to 0$
and $\beta _n \sin \zeta _n$ converges (to $0$ or $t$ 
depending on if $s_n ^2 +a_n ^2 \to \infty$ or not).  Since $\beta _n \to \infty$, we also get that 
$(\frac{s_n^2}{a_n ^2}+a_n ^2) \sin \zeta _n  \to 0$.  Since $\zeta _n \to 0$ and $\beta _n \to \infty$, we get that $\mathfrak{s}_n \mathfrak{a}_n Z_n \mathfrak{a}_n ^{-1} \mathfrak{s}_n ^{-1} \to I$ for $Z_n = \begin{pmatrix}
    \cos \zeta_n & \sin \zeta_n\\
    -\sin \zeta_n & \cos \zeta_n
\end{pmatrix} \in K$.  Since $\zeta _n \to 0$, $( 1-\cos \zeta _n ) / \sin \zeta _n \to 0$. Hence  
\[\alpha _n (1-\cos \zeta _n)= (\beta _n \sin \zeta _n)(\frac{\alpha _n }{\beta _n})(\frac{1-\cos \zeta _n }{\sin \zeta _n}) \to 0\] and 
\begin{align*}
-\alpha _n \frac{\sin \zeta _n}{a_n ^2} + \beta _n \big(1-\cos \zeta _n -\frac{s_n}{a_n ^2}\sin \zeta _n \big) &\\=  \frac{-\alpha _n} {\beta _n}\frac{\beta _n \sin \zeta _n}{a_n^2} 
+\beta _n \sin \zeta _n (\frac{1-\cos \zeta _n}{\sin \zeta _n })- \beta _n \frac{s_n}{a_n ^2} \sin \zeta _n &\\\to 0.
\end{align*}
Thus, 
\[
(\mathfrak{s}_n\mathfrak{a}_n,v_n)\ Z_n\ (\mathfrak{s}_n\mathfrak{a}_n,v_n)^{-1} \to (I, (t,0)).
\]

Now, replacing $\Theta _n $ by $\Theta _n Z _n^{-1}$, we may assume that $x_n \to 0$.   
Since $(\frac{s_n^2}{a_n ^2}+a_n ^2) \sin \theta _n \to 0$ and $\frac{s_n ^2}{a_n ^2}+a_n ^2 
\not \to 0$, we get that $\sin \theta _n \to 0$.  Replacing $\theta _n$ by a translate of $2k\pi$, we may assume that $\theta _n\to 0$.  This implies that $\frac{1-\cos \theta_n}{\sin \theta _n }\to 0$.  Since $x_n \to 0$, 
we get that $\beta _n (\frac{s_n ^2}{a_n^2 }+ a_n ^2) \sin \theta _n \to 0$, hence 
$\alpha _n (1-\cos \theta _n +\frac{s_n}{a_n ^2}\sin \theta _n ) \to 0$.  
Therefore, $$|\alpha _n \frac{- \sin \theta _n}{a_n^2} | \leq |\frac{\alpha _n}{\beta _n}| |
\beta _n (\frac{s_n ^2}{a_n ^2}+a_n ^2) \sin \theta _n  | \to 0$$ and 
$$\begin{array}{cl}  & |\beta _n (1-\cos \theta _n - \frac{s_n}{a_n ^2}\sin \theta _n ) |   \\ = &
|\frac{1-\cos \theta _n}{\sin \theta _n }\frac{a_n ^2}{s_n ^2 +a_n ^4} - \frac{s_n}{s_n ^2+a_n ^4} |
|\beta _n (\frac{s_n ^2}{a_n ^2}+a_n ^2) \sin \theta _n | \to 0. \end{array}  $$
This implies that $H= V_0= \{ (I, (x, 0)) \mid x\in \R\}$.

\noindent {\bf Case $2$: } Suppose $\frac{\beta _n}{\alpha _n} \to 0$.  
Let $F= \begin{pmatrix} 0 & -1 \\ 1& 0 \end{pmatrix}$.  Since $(g_n , v_n) K (g_n, v_n)^{-1} \to H$, 
conjugating by $F$, we get that $(Fg_n, Fv_n) K (Fg_n, Fv_n)^{-1} \to FHF^{-1}$.  Since $Fv$ satisfies case (i) we get that $FHF^{-1} = V_0$.  Thus, $H= V_\infty = \{ (I, (0,y) \mid y\in \R \}$.

\noindent {\bf Case $3$: } Suppose $c\not =0$ is a limit point of $\frac{\alpha _n}{\beta _n}$.  
By passing to a subsequence we may assume that $\frac{\alpha _n}{\beta _n} \to c$.  
There are two cases here:  if one of $a_n $ or $s_n$ diverges to $\infty$, then as in case (i) we can show 
that $H= V_0$.   We now assume that $a_n =1$ and $s_n=0$.  Then, choosing $\zeta _n\rightarrow 0$ such that $\beta_n\sin \zeta_n \rightarrow t$ for any preassigned real number $t$, it is easy to see that 
\[g_nZ_n g_n^{-1}\rightarrow \big(I,\ (t,-ct)\big). \] for suitable $Z_n \in K$.    
Therefore, following a similar argument as in case (i), we may show in this case that
$H= V_c=\{ (I, (t, ct)) \mid t \in \R \}$.

We now consider the case when $(v_n)$ has a convergent subsequence.  In this case, by passing to a 
subsequence and conjugating by $(I, -v_n)$, we may assume that $(\mathfrak{s}_n\mathfrak{a}_n, 0) K(\mathfrak{s}_n\mathfrak{a}_n ,0)^{-1} \to H$.  
Since $\mathfrak{s}_n , \mathfrak{a} _n \in SL(2, \R)$ and $K\subset SL(2, \R)$, $H\subset SL(2, \R)$.  
Recall that we also have $s_n \to \infty$ or $s_n=0$ and $a_n \to \infty$ or $a_n=1$.  
If both $s_n=0$ and $a_n =1$, we are done.  So we may assume that either $s_n \to \infty $ or $a_n \to \infty$.  
We consider the case that $s_n \to \infty$ and other case that $a_n \to \infty$ is similar.  
 We have 
 \[
 T_n = \begin{pmatrix} \cos \theta _n - \frac{s_n}{a_n ^2}\sin \theta _n  & (\frac{s_n ^2}{a_n^2}+a_n ^2) \sin \theta _n \\ {-\sin \theta _n / a_n ^2} & \cos \theta _n +\frac{s_n}{a_n ^2}\sin \theta _n \end{pmatrix} \to g = \begin{pmatrix} g_{1,1} & g_{1,2} \\ g_{2,1} & g_{2,2} 
\end{pmatrix}, {~~\rm say}.
\]
This, in particular, implies that $(\frac{s_n^2}{a_n ^2}+a_n ^2) \sin \theta _n$ 
converges.  Since $\frac{s_n}{a_n^2} +a_n^2 \to \infty$, $\sin \theta _n \to 0$ and since 
$\frac{\sin \theta _n}{a_n^2 }\to g_{2,1}$, we get that $g_{2,1}=0$.  
Also, for $t\in \R$, we choose 
$\zeta _n\to 0$ such that $(\frac{s_n^2}{a_n ^2}+a_n ^2)\sin \zeta _n \to t$.  Then $\cos \zeta _n\to 1 $ and 
$\frac{s_n}{a_n ^2}\sin \zeta _n \to 0$, hence $g_n Z_ng_n^{-1} \to \begin{pmatrix} 1  & t \\ 
0 & 1\end{pmatrix}  $ for suitable $Z_n\in K$. If, in addition $g_{1,2}=0$, then
$(\frac{s_n ^2}{a_n^2}+a_n ^2) \sin \theta _n \to 0$ and hence $\frac{s_n}{a_n ^2}\sin \theta _n \to 0$.  
Thus, $\cos \theta _n$ converges to $g_{1,1}$ as well as to $g_{2,2}$. But, $g_{1,1} g_{2,2}=1$, therefore, 
$g_{1,1}=\pm 1 = g_{2,2}$.  In fact, taking $\theta _n = \pi$, 
we get that 
 \[ 
g_n\Theta_ng_n^{-1}=\begin{pmatrix}
    \cos \theta_n - \frac{s_n}{a_n^2}\ \sin \theta_n & (\frac{s_n^2}{a_n^2}+a_n^2)\  \sin \theta_n\\
    -\sin \theta_n /a_n^2 & \cos \theta_n + \frac{s_n}{a_n^2} \sin \theta_n
\end{pmatrix} = \begin{pmatrix} -1 &0 \\ 0&-1 \end{pmatrix}. 
\]
Thus, we have $H= \tilde {N^+}$, that is $g_nKg_n^{-1}\rightarrow \tilde {N^+}$.\\\\
\end{proof}
\section{Conjugacy limits of maximal diagonalizable subgroups}
Now we study limit points of the set of maximal diagonalizable subgroups of $G=\SL(2,\R)\ltimes \R^2$. Recall that $A= \left\{ \begin{pmatrix} a &0 \\ 0 &1/a \end{pmatrix} \mid a\in \R^+ \right\}$.
We consider the conjugacy limits $(g_n,v_n)A(g_n,v_n)^{-1}$ of $A$ for different types of $g_n$ and $v_n$. 
\begin{proof}[Proof of Theorem \ref{limit-D}] 
Let $H$ be a limit point of 
$\{ (g, v) A (g,v)^{-1} \mid g\in SL(2, \R), ~ v\in \R ^2 \}$. So, there are $g_n \in SL(2, \R)$ and $v_n \in \R ^2$ such that $(g_n, v_n) A (g_n,v_n)^{-1} \to H$. Using Iwasawa decomposition, we may assume that 
$g_n\in N^+$.  Let $g_n = \begin{pmatrix} 1 & s_n \\ 0 & 1 \end{pmatrix}$ and $v_n = (\alpha _n , \beta _n)$.  
Suppose $(g, v) \in H$. Then, there are $a_n \in \R ^+$ such that $(g_n, v_n) \mathfrak{a}_n(g_n,v_n)^{-1} \to (g, v)$, where $\mathfrak{a}_n = \begin{pmatrix} a_n &0 \\ 0 &1/a_n \end{pmatrix}$.  This implies that 
\[
\begin{pmatrix} a_n & ( \frac{1}{a_n}-a_n)s_n \\ 0 & \frac{1}{a_n} \end{pmatrix} \to g, ~~{\rm and} ~~
\]
\[
((1-a_n ) \alpha _n + (a_n - \frac{1}{a_n})s_n \beta _n, (1-\frac{1}{a_n})\beta _n ) \to v.
\]
By passing to subsequences and taking conjugates, we may assume that $s_n=0$ or $s_n \to \infty$, $\alpha _n=0$ or 
$\alpha _n \to \infty$ and $\beta _n=0$ or $\beta _n \to \infty$.  

We first consider the case that $s_n\to \infty$.  
If both $\alpha _n=0$ and $\beta _n=0$, then it is easy to see that 
$H = N^+=\left\{ \begin{pmatrix} 1 & t \\ 0&1 \end{pmatrix} \mid t\in \R \right\}$.  So assume that either $\alpha _n \to \infty$ or $\beta _n \to \infty$.  
Then $v_n\to \infty$.  This implies by Lemma \ref{el1} that $g$ is unipotent.  Therefore, $a_n \to 1$.

Since $g_n \in N^+$, and $B =AN^+$, it follows that $(g_n, v_n)A(g_n,v_n)^{-1} $ are closed subgroups of $B\ltimes \R ^2$.  Since $B\ltimes \R ^2$ is isomorphic to a group of upper triangular matrices with positive diagonal entries, $B\ltimes \R ^2$ is an exponential Lie group. Hence, by Proposition \ref{cs}, $H$ is a connected one-dimensional subgroup of $N^+\ltimes \R ^2$ which is the $3$-dimensional Heisenberg group.  
Therefore, 
\[
H= V_{(a,b,c)}= \left\{ \left ( \begin{pmatrix} 1 & at \\ 0 & 1 \end{pmatrix} ,\left( ct+\frac{abt^2}{2}, bt \right) \right) \mid t\in \R \right\} 
\]
for some 
$a, b, c \in \R$.  Then for each $t\in \R$, there are $b_n$ (depending on $t$) such that 
$$b_n \to 1, ~~ (\frac{1}{b_n} -b_n) s_n \to at ~~{\rm and}~~ (1- \frac{1}{b_n})\beta _n \to bt .$$

\noindent {\bf Case $1$:} Suppose $\frac{s_n}{\beta _n} \to r \in \R$.  Then $a=-2br$ and hence 
$$(1-b_n)\alpha _n - \frac{(1+b_n)}{b_n} (1-b_n)s_n \beta _n \to ct-rb^2t^2 .$$  If $b=0$, then $H= V_0$.  
So, we may assume that $b\not = 0$.  
Also $$2(1-b_n)s_n \beta _n  - \frac{(1+b_n)}{b_n} (1-b_n )s_n \beta _n \to -rb^2 t^2.$$  Therefore, 
$$(1-b_n) \alpha _n -2(1-b_n)s_n \beta _n \to ct .$$  This implies that $2s_n - \frac{\alpha _n}{\beta _n} \to \frac{c}{b}$.  
Let $d= \lim 2s_n - \frac{\alpha _n}{\beta _n}$.  Then $c= db$.  Thus, $H= V_{(r_1b, b, r_2b)}$ for some $r_1, r_2 \in \R$. 
Note that if $b\not =0 $, $V_{(r_1b, b, r_2b)}= V_{(a, 1, c)}$ for some $a, c\in \R$. 

\noindent{\bf Case $2$:} Suppose $\frac{\beta _n}{s_n }\to p \in \R$,  Then $b=- \frac{pa}{2}$ and hence 
$$(1-b_n)\alpha _n - \frac{(1+b_n)}{b_n} (1-b_n)s_n \beta _n \to ct- \frac{pa^2t^2}{4} .$$  
If $a=0$, then $H = V_0$.  So, we may assume that $a\not =0$.  Thus, $$(1-b_n) \alpha _n - 2(1-b_n) s_n \beta _n \to ct.$$  
This implies that $\frac{\alpha _n}{s_n} -2\beta _n \to 2\frac{c}{a}$, hence $c= \frac{da}{2}$ where 
$d= \lim \frac{\alpha _n}{s_n} - 2 \beta _n$.  Thus, 
$H= V_{(a, r_1a, r_2a)}$ for some $r_1, r_2 \in \R$.  Note that if $a \not =0$, then $V_{(a, r_1a, r_2a)}= V_{(1, b,c)}$ for some 
$b, c\in \R$.

We now consider the case $s_n=0$.  Then $g= I$ and $$((1-a_n ) \alpha _n, (1- \frac{1}{a_n})\beta _n ) \to v$$ 
for some $a_n \to 1$.

\noindent{\bf Case $3$:} Suppose $\frac{\beta _n}{\alpha _n} \to 0$.  Then since $v_n = (\alpha _n , \beta _n ) \to \infty$, 
$\alpha _n \to \infty$ and 
$(1- \frac{1}{a_n})\beta _n = \frac{-1}{a_n}(1-a_n)\alpha _n \frac{\beta _n}{\alpha _n}\to 0$.  Therefore 
$H\subset V_0$.  For any $t\in \R$, taking $b_n = 1- \frac{t}{\alpha _n}$, we get that $b_n \to 1$, $(1-b_n)\alpha _n =t$.  
Therefore, $((1-b_n ) \alpha _n, (1- \frac{1}{b_n})\beta _n )\to (t,0)$.  Thus, $H=V_0 = V_{(0,1,1)}$.  

\noindent{\bf Case $4$:} Suppose $\frac{\alpha _n}{\beta _n} \to 0$.  Then 
$(1-a_n)\alpha _n = -a_n(1- \frac{1}{a_n})\beta _n \frac{\alpha _n}{\beta _n}\to 0$.  Therefore, 
$H\subset V_\infty = \{ (0, t) \mid t\in \R \}$.  For any $t\in \R$, take 
$b_n = \frac{\beta _n}{\beta _n -t}$.  Then $b_n \to 1$ and $((1-b_n ) \alpha _n, (1- \frac{1}{b_n})\beta _n )\to (0,t)$.  Thus, $H= V_\infty = V_{(0,1,0)}$.

\noindent{\bf Case $5$:} Suppose $\frac{\beta _n}{\alpha _n }\to c \not = 0$.  
Let $x,y \in \R$ be such that $v=(x,y)$.  Then 
$(1-a_n ) \alpha _n \to x$
and $(1- \frac{1}{a_n})\beta _n \to y$.  
Since $(1- \frac{1}{a_n})\beta _n = \frac{-1}{a_n} (1-a_n ) \alpha _n  \frac{\beta _n}{\alpha _n}$, we get that $y=-cx$.  Therefore $H\subset V_r= \{ (t, rt ) \mid t\in \R \}$ for $r=-c$.  

For $t\in \R$, let $b_n = 1- \frac{t}{\alpha _n}$.  Then $b_n \to 1$ and $(1-b_n)\alpha _n =t$.  
Also, $(1- \frac{1}{b_n})\beta _n = \frac{-t}{\alpha _n -t }\beta _n \to -tc$.  Therefore 
$V_r \subset H$.  Thus, $H= V_r = V_{(0,1,r)}$.
\end{proof}
\section{Conjugacy limits of unipotent and Borel subgroups}

We now find the conjugacy limits of the Borel subgroup 
\[
B:=\left\{
\left(\begin{pmatrix} a&s\\ 0&1/a\end{pmatrix}:\ a,s\in\R\right)
\right\}
\]
of $\SL(2,\R)$ inside $G=\SL(2,\R)\ltimes \R^2$.
\begin{proof}[Proof of Theorem \ref{limits-B}]
Clearly, a conjugacy limit of $B$, when conjugated by elements of the type $(g,0)$ with $g\in\SL(2,\R)$, is $B$ itself or a conjugate of $B$ by elements of $K$. For conjugation by more general type of elements, it is enough to consider conjugation by elements of the form $(I,v)$ with $v\in\R^2$. So, let   $\mathfrak{g}_n=(I,v_n)$ and $\mathfrak{b}_n=\begin{pmatrix}a_n&s_n\\0&1/a_n\end{pmatrix}$. Then a simple calculation shows that
\begin{align*}
\mathfrak{g}_n\mathfrak{b}_n\mathfrak{g}_n^{-1}=\left(\begin{pmatrix}a_n&s_n\\0&1/a_n\end{pmatrix},\Big((1-a_n)\alpha_n -s_n\beta_n ,\ (1-1/a_n)\beta_n\Big)\right)\\ =\big(\mathfrak{b}_n,(x_n,y_n)\big) ,\  (\text{say})
\end{align*} where $v_n = (\alpha _n , \beta _n)$.
Assume that $v_n \to \infty$ and $\mathfrak{g}_n\mathfrak{b}_n\mathfrak{g}_n^{-1}\to (g, v)$.  Then by Lemma \ref{el1}, $g$ is unipotent, hence $a_n\rightarrow 1$ necessarily. 
We shall consider the following cases.

\noindent {\bf Case $1$:} Suppose $\frac{\alpha_n}{\beta_n}\rightarrow 0$. Then $\beta _n \to \infty$.  
If $s_n\rightarrow s$ for some non-zero $s$, then it is easy to see that $x_n$ cannot converge. Hence, $s_n\rightarrow 0$. Given $x,y\in \R$, choose $s_n\rightarrow 0$ and $a_n\rightarrow 1$ in such a way that $-\beta_n s_n \rightarrow x$ and $(1-1/a_n)\beta_n\rightarrow y$.  Then 
$(1-a_n)\alpha _n = -a_n (1- \frac{1}{a_n})\beta _n \frac{\alpha _n}{\beta _n} \to 0$, hence 
$\mathfrak{g}_n \mathfrak{b}_n ' \mathfrak{g}_n^{-1}\rightarrow \big(I,(x,y) \big)$ for suitable $\mathfrak{b}_n '\in B$.  This shows that, in this case, the limit subgroup is $\R^2$.  

\noindent{\bf Case $2$:} Suppose $\frac{\beta_n}{\alpha_n}\rightarrow  0$. Then $\alpha _n \to \infty$.  
Let $x\in \R$.  Then for $\tilde {\mathfrak{b} }_n = \begin{pmatrix} 1- \frac{x}{\alpha _n} & 0 \\ 0 & \frac{\alpha _n}{\alpha _n -x} \end{pmatrix}$, we get that $\mathfrak{g}_n\tilde {\mathfrak{b}_n} \mathfrak{g}_n^{-1} \to (I, (t, 0))$.  
Thus, $\{ (I, (t,0))\mid t \in \R \} \subset H$.  

If $\beta _n$ converges, then $y_n =(1- \frac{1}{a_n})\beta _n\to 0$ and it can be shown that $H=N^+\ltimes \{(t, 0) \mid t\in \R \}$.  So, we may assume that $\beta _n \to \infty$.  Since $\mathfrak{b}_n$ converges, $s_n$ converges, say to $s$.  Since $\beta _n \to \infty$ and $x_n$ converges, we get that $(1-a_n) \frac{\alpha _n}{\beta _n} \to  s $.

If $\frac{\beta _n ^2}{\alpha _n} \to c \in \R$, then $y_n =\frac{-1}{a_n} (1-a_n) \frac{\alpha _n}{\beta _n } \frac{\beta _n ^2}{\alpha _n} \to -sc$.  Therefore, by Proposition \ref{cs} we get that 
$H= \{ (\mathfrak{s}, (t,-cs) \mid t, s \in \R \}$ where $\mathfrak{s} = \begin{pmatrix} 1& s \\ 0 & 1 \end{pmatrix}$.  

If $\frac{\beta _n ^2}{\alpha _n}\to \infty$, then $(1-a_n) \frac{\alpha _n}{\beta _n} = (1- \frac{1}{a_n})\beta _n 
(-a_n) \frac{\alpha _n}{\beta _n ^2}\to 0$.  Therefore, by Proposition \ref{cs} we get that 
$H= \R ^2 = \{ (I, (t,r) \mid t, r \in \R \}$. 

\noindent{\bf Case $3$}: Suppose $\frac{\alpha_n}{\beta_n}
\rightarrow c$ for some non-zero real number $c$. Note that we are necessarily assuming both $|\alpha_n|$ and $|\beta_n|$ are unbounded. Otherwise, the limit subgroup will be a conjugate of $B$. Now, suppose $x_n\rightarrow x$. Then $\frac{x_n}{\beta_n}\rightarrow 0$. But, 
\[
\frac{x_n}{\beta_n}=\frac{(1-a_n)\alpha_n}{\beta_n}-s_n
\]
which implies that $s_n\rightarrow 0$ in this case as well.  Therefore $H\subset \R ^2$.  Thus, by Proposition \ref{cs}, $H$ is of dimension $2$ and hence $H= \R ^2$.  
\end{proof}
We next determine the conjugacy limits of the unipotent subgroup $N^+$.
\begin{proof}{Proof of Theorem \ref{limits-U}}
Since $N^+$ is normal in $B$, using Iwasawa decomposition in $SL(2,\R)$, it follows that whenever $g_n\in\SL(2,\R)$, the limit of $g_nN^+g_n^{-1}$ is either $N^+$ or its conjugates by elements of $K$. For conjugation by other elements of $G$, without loss of generality, we may consider conjugation by elements of the form $\mathfrak{g}_n=(k_n,v_n)$ with $k_n\in K$ and $v_n=(\alpha_n,\beta_n)\in\R^2$. Also, it is easy to see that if $L$ is a closed subgroup of $G$, $g_nLg_n^{-1}\rightarrow H$ and $k_n\rightarrow k$, then $k_ng_nLg_n^{-1}k_n^{-1}\rightarrow kHk^{-1}$. In view of this, it is enough to consider conjugation of $N^+$ by elements of the form $\mathfrak{g}_n=(I,v_n)$. 

Now, suppose $\mathfrak{g}_n=\big(I,v_n\big)$ with $v_n=(\alpha_n,\beta_n)$. Then for any $\mathfrak{s}_n=\begin{pmatrix}
    1&s_n\\0&1
\end{pmatrix}$, 
\[
\mathfrak{g}_n\mathfrak{s}_n\mathfrak{g}_n^{-1}=\left(\begin{pmatrix}
    1&s_n\\0&1
\end{pmatrix},\big(-s_n\beta_n,0\big)\right).
\]
 If $|\beta_n|\rightarrow\infty$, then choosing $s_n\rightarrow 0$ suitably $\big(I,(t,0)\big)$ can be realized as a limit of $\mathfrak{g}_n\mathfrak{s}_n\mathfrak{g}_n^{-1}$ for any real number $t$.  Also, if $\mathfrak{g}_n \mathfrak{s}_n\mathfrak{g}_n^{-1}$ converges, then 
 $s_n$ as well as $s_n\beta _n$ converges. but, as $|\beta _n| \to \infty$, we have $s_n \to 0$.  This shows that, in this case, the limit of $\mathfrak{g}_nN^+\mathfrak{g}_n^{-1}$ is $V_0$. 
\end{proof}
\acknowledgement{
Both authors thank the International Centre for Theoretical Sciences, Bangalore, for its hospitality during the programme "Ergodic Theory and Dynamical Systems" from $5$th December to $16$th December $2022$, when this work was initiated. The first author thanks the Indian Statistical Institute, Bangalore Centre, for its hospitality during his numerous visits while this work was in progress. The second author acknowledges the support from SERB through the grant under MATRICS MTR/2022/000429.}

\bibliography{refnil}
\bibliographystyle{plain}

\end{document}